\input amstex
\documentstyle{amsppt}
\NoRunningHeads
\magnification=\magstep1
\baselineskip=12pt
\parskip=5pt
\parindent=18pt
\topskip=10pt
\leftskip=0pt
\pagewidth{32pc}
\pageheight{47pc}
\topmatter
\title Projectivity via the dual K\"ahler cone - Huybrechts' criterion 
\endtitle
\author Keiji Oguiso and Thomas Peternell 
\endauthor

\abstract In this note we give an elementary proof 
for a remarkable criterion due to Daniel Huybrechts for a K\"ahler 
surface to be projective. 
\endabstract

\endtopmatter

\document
\head Introduction \endhead 
One of the main idea in higher dimensional algebraic geometry is to 
study varieties through numerical properties of their cones, of which 
origin probably goes back to the Kleiman criterion, 
the duality between the ample cone and the so-called Kleiman-Mori cone, 
the cone of effective curves (cf. [KMM]). 
\par
\vskip 4pt
Quite recently, Daniel Huybrechts took this idea into his study of 
hyperk\"ahler manifolds and stated as a byproduct 
the following remarkable criterion to distinguish projective 
surfaces from the cone theoretical view point: 

\proclaim{Huybrechts' Criterion ([Hu, Remark 3.12 (iii)])} 
A compact K\"ahler surface is projective if and only if the dual 
cone of the K\"ahler cone contains an inner integral point. 
(For the precise definitions, see (1.4) in section 1.)  
\endproclaim   

However his original proof relies on powerful 
but highly advanced techniques in complex analysis (Demailly's 
singular Morse theory) and he himself asked in the same paper 
whether it is possible or not to prove this in a more elementary way.
\par
\vskip 4pt
The aim of this short note is to answer his question by 
giving a proof based on more or less familiar results 
on surfaces found now in standard books, [Beu], [BPV] and [GH]. 
Our proof is based on the notion of algebraic dimension 
while it is almost free from the classification of surfaces.  
\par \vskip 4pt Of course it is very interesting to ask whether Huybrechts' criterion
also holds in higher dimensions. We will address to this question in a second
part of this paper [OP]. 

\head Acknowledgement \endhead
The first named author would like to express his sincere thanks to Professors 
H\'el\`ene Esnault and Eckart Viehweg for their invitation 
to Germany, Professor Yujiro Kawamata for his encouragement 
and enlighting suggestion during his stay in Moscow June 1999. 
Both authors would like to express their thanks to Daniel Huybrechts 
for his interest in this work.
This work has been done during the first named author's stay in Universit\"at 
Bayreuth July 1999 under the financial support by the Alexander-Humboldt 
fellowship and the Universit\"at Bayreuth. The first named author would like 
to express his best thanks for their financial support which made this work 
possible.  

\head 1.Preliminaries \endhead
\demo{(1.0)} Thoughout this note, the term surface means a compact, connected 
complex manifold of dimension two. Let $S$ be a surface. 
The transcendental degree of the meromorphic function field of $S$ over 
$\Bbb C$ is called the 
algebraic dimension and is denoted by $a(S)$. It is well known that 
$a(S) \in \{0, 1, 2\}$ and $S$ is projective if and only if $a(S) = 2$.
\enddemo
\demo{(1.1)} A Hermitian metric $g$ on $S$ is called K\"ahler 
if the associated positive real $(1,1)$ form $\omega_{g}$ 
is $d-$closed. We call $\omega_{g}$ a K\"ahler form if $g$ is 
a K\"ahler metric. A surface is called K\"ahler if it admits at least 
one K\"ahler metric. Note that every projective surface is K\"ahler 
but the converse is not true in general. \enddemo
Let $S$ be a K\"ahler surface.
\demo{(1.2)} By definition, any K\"ahler metric $g$ on $S$ determines 
a de Rham cohomology class $[\omega_{g}]$. This class lies in 
the real $(1,1)$ part $H^{1,1}(S, \Bbb R)$ of the Hodge 
decomposition of $H^{2}(S, \Bbb C)$. We often abbreviate 
$H^{1,1}(S, \Bbb R)$ by $H^{1,1}$.
We call an element $\eta \in H^{1, 1}$ a K\"ahler class 
if it is represented by a K\"ahler form, that is, in the case where 
there exists a  K\"ahler metric $g$ such that $\eta = [\omega_{g}]$. 
\enddemo   
\demo{(1.3)} The real vector space $H^{1, 1}$ carries a natural symmetric 
bilinear form $(*.*)$ induced by the cup product on 
the integral cohomology group $H^{2}(S, \Bbb Z)$. It is well known that 
$(*.*)$ on $H^{1, 1}$ is non-degenerate and is of 
signature $(1, h^{1,1}(S)-1)$. We also regard the finite dimensional 
real vector space $H^{1, 1}$ as a linear topological space 
using some norm $\vert * \vert$. Therefore we can speak of 
the closure $\overline{A}$ of $A \subset H^{1, 1}$. For $a \in H^{1, 1}$ 
and for a positive real number $\epsilon 
> 0$, we set 
$$B_{\epsilon}(a) := \{ x \in H^{1,1} \vert \vert x - a \vert 
\leq \epsilon\}.$$ 
Furthermore, let $U_{\epsilon}(a)$ be the interior of $B_{\epsilon}(a)$ 
and $\partial B_{\epsilon}(a)$ its boundary.  
\enddemo
\demo{(1.4)}
\roster
\item 
The K\"ahler cone $\Cal K(S)$ of $S$ is the subset of $H^{1,1}$ 
consisting of the K\"ahler classes of $S$. 
By definition, $\Cal K(S)$ is a convex cone of $H^{1,1}$. 
It is also well known that $\Cal K(S)$ is an open subset of $H^{1,1}$.  
\item 
The dual cone $\Cal K^{*}(S)$ of the K\"ahler cone $\Cal K(S)$ 
is the set of elements $x \in H^{1,1}$ such that $(x.\eta) > 0$ 
for any $\eta \in \Cal K(S)$. 
\item 
An element $x$ of $\Cal K^{*}(S)$ is called integral if 
$x \in \Cal K^{*}(S) \cap \iota^{*}H^{2}(S, \Bbb Z)$, where 
$\iota : \Bbb Z \rightarrow \Bbb R$ is a natural inclusion of sheaves. 
An integral element is nothing but an element of $\Cal K^{*}(S) \cap NS(S)$ 
(cf.(1.7)).
\item 
An element $x$ of $\Cal K^{*}(S)$ is called an inner point if there 
exists a positive real number $\epsilon > 0$ such that 
$U_{\epsilon}(x) \subset \Cal K(S)^{*}$.
\endroster
\enddemo
\proclaim{Lemma (1.5)} Let $(H, |*|)$ be a finite dimensional, real normed 
vector space equipped with a real valued, 
non-degenerate bilinear form $(*.*)$. 
Let $K \subset H$ be a non-empty convex subset such that $0 \not\in K$. 
Set $K^{*} \subset H$ to be the dual of $K$ with respect to $(*.*)$. 
Let $x \in H$. Then $x$ is an inner point of $K^{*}$ if and only if 
there exists a positive real number $r > 0$ such that 
$(x, \eta) \geq r \vert \eta \vert$ for all $\eta \in \overline{K}$. 
\endproclaim 
\demo{Proof} This will follow from the compactness of the space 
$B_{\epsilon}(x) \times (\overline{K} \cap \partial B_{1}(0))$. \qed \enddemo  
The following direct consequence is crucial for our proof: 

\proclaim{Corollary (1.6)} 
Let $x \in H^{1,1}$. Then $x$ is an inner point of $\Cal K^{*}(S)$ 
if and only if there exists a positive real number $r > 0$ such that 
$(x. \eta) \geq r \vert \eta \vert$ for all $\eta 
\in \overline{\Cal K(S)}$. 
\qed
\endproclaim
         
\demo{(1.7)} The group $H^{1,1} \cap \iota^{*}H^{2}(S, \Bbb Z)$ is called 
the N\'eron-Severi group of $S$ 
and is denoted by $NS(S)$. The rank of $NS(S)$ is called the Picard number of 
$S$ and is written by $\rho(S)$. By the Lefschetz $(1,1)$ Theorem, 
each element of $NS(S)$ is represented by the first Chern 
class of some line bundle. However, contrary to the projective case, 
the natural map from the 
group of Cartier divisors to the Picard group 
is not surjective in general. So, we CAN NOT say that each element 
of $NS(S)$ is represented by a divisor in the K\"ahler category.   
\enddemo

\head 2. K\"ahler cones of K3 surfaces and complex tori \endhead
\proclaim{Theorem (2.1) [Bea, Page 123, Theorem 2]} 
Let $S$ be a K3 surface, that is, a (K\"ahler) surface such that 
$K_{S} = 0$ in $\text{Pic}(S)$ and that $\pi_{1}(S) = \{1\}$. Let 
$\Cal C^{+}(S)$ be the connected component of the space 
$\{x \in H^{1,1} \vert (x.x) > 0\}$ which contains the K\"ahler 
classes. Then the K\"ahler cone $\Cal K(S)$ coincides with 
the subspace $\widetilde{\Cal K(S)}$ of $\Cal C^{+}(S)$ defined by 
$(x.[C]) > 0$ 
for all non-singular rational curves $C$ in $S$, that is, 
\par 
$\Cal K(S) = \widetilde{\Cal K(S)} := \{x \in \Cal C^{+}(S) 
\vert (x.[C]) > 0$ 
for all $C \simeq \Bbb P^{1}$ in $S$ $\}$. \qed 
\endproclaim 
\demo{Remark} It is clear that $\Cal K(S) \subset \widetilde{\Cal K(S)}$. 
However, the other inclusion $\widetilde{\Cal K(S)} \subset \Cal K(S)$
is highly non-trivial. For details, we refer to [Bea]. 
\enddemo
\proclaim{Theorem (2.2)}   
Let $S$ be a complex torus of dimension 2. Let 
$\Cal C^{+}(S)$ be the connected component of the space 
$\{x \in H^{1,1} \vert (x.x) > 0\}$ which contains the K\"ahler 
classes. Then $\Cal K(S) = \Cal C^{+}(S)$. 
\endproclaim
\demo{Remark} This result should be known. However, the authors 
could not find any references. The present proof was communicated to us by
D. Huybrechts; our original proof is more complicated and works by reduction to
the algebraic case. 
\enddemo
\demo{Proof} First notice that any $(1,1)-$ class can be represented by a form
with constant coefficients. Suppose $\Cal K(S)  \not = \Cal C^{+}(S).$
Since $\Cal K(S) \subset  \Cal C^{+}(S),$ we find a constant $(1,1)-$form $\phi$
such that $[\phiÊ] \in \Cal C^{+}(S) \cap \partial \Cal K(S).$ Then $\phi$ is semipositive
but not positive. Therefore $\phi^2 = 0,$ contradiction. \qed
\enddemo

\proclaim{Lemma (2.3)} Let $S$ be a minimal K\"ahler surface. 
Assume that $a(S) = 0$. Then $S$ is either a K3 surface or a 
complex torus of dimension 2. 
\endproclaim 
\demo{Proof} This is of course well known, see e.g. [BPV]. 
We give a proof to convince the reader that no deep result from 
classification theory is involved. Since $a(S) = 0$, 
we have $\kappa(S) = 0$ or $-\infty$, where $\kappa(S)$ 
is the Kodaira dimension of $S$. Moreover, if $h^{0}(K_{S}) = 0$, then 
by the Serre duality $h^{2}(\Cal O_{S}) = 0$ and $S$ is 
then projective by the Kodaira criterion, 
a contradiction. Therefore $K_{S} = 0$ in $\text{Pic}(S)$ by the minimality 
of $S$. Since $S$ is K\"ahler, this gives the result. \qed
\enddemo
In order to prove Huybrechts' criterion, we also need to know the structure 
of the N\'eron-Severi groups of K3 surfaces and complex tori of algebraic 
dimension zero.         
  
\proclaim{Proposition (2.4)} Let $S$ be a K3 surface. Assume that 
$a(S) = 0$. Then, 
\roster
\item $\text{Pic}(S)$ and $NS(S)$ are torsion free 
and are isomorphic under $c_{1}$. Moreover 
$NS(S) \otimes \Bbb R$ is negative definite with respect to 
$(*.*)$. 
\item $S$ contains at most 19 distinct smooth rational curves and 
contains no other curves. 
\endroster
\endproclaim 
\demo{Proof of (1)} The first part of the assertion is well known. 
Using $a(S) = 0$ and the Riemann-Roch Theorem, we readily see that 
$L^{2} < 0$ for all $L \in \text{Pic}(S) - \{0\}$. 
Since $(*.*)$ is defined over $\Bbb Z$, this implies the result. \qed \enddemo

\demo{Proof of (2)} Let $C$ be an irreducible curve on $S$. Then 
$C \simeq \Bbb P^{1}$, because $0 > C^{2} = (K_{S} + C. C) 
= 2p_{a}(C) - 2$ by (1) and the adjunction formula.
\proclaim{Claim 1} Let $C_{1}, ..., C_{m}$ be $m$ distinct irreducible 
curves on $S$. Then $[C_{1}], ..., [C_{m}]$ are 
linearly independent in $NS(S) \otimes \Bbb R$. \endproclaim 
\demo{Proof} Since the classes $[C_{i}]$ defined over 
$\Bbb Z$, it is enough to show that if $\sum_{i \in I} p_{i}[C_{i}] = 
\sum_{j \in J} q_{j}[C_{j}]$, where $I \cap J = \emptyset$, 
$p_{i} \in \Bbb Z_{\geq 0}$ and $q_{j} \in \Bbb Z_{\geq 0}$ 
then $p_{i} = q_{j} = 0$. 
\par 
By (1), we have 
$$0 \geq (\sum_{i \in I} p_{i}[C_{i}].\sum_{i \in I} p_{i}[C_{i}]) 
=  (\sum_{i \in I} p_{i}[C_{i}].\sum_{j \in J} q_{j}[C_{j}]) \geq 0.$$
 Therefore, 
$(\sum_{i \in I} p_{i}[C_{i}].\sum_{i \in I} p_{i}[C_{i}]) = 0$. 
Then again by (1), we have 
$$[\sum_{i \in I} p_{i}C_{i}] = \sum_{i \in I} p_{i}[C_{i}] = 0$$ 
in $NS(S)$ and $\sum_{i \in I} p_{i}C_{i} = 0$ in $\text{Pic}(S)$. 
This is possible only in the case where $p_{i} = 0$ for all $i \in I$. 
Similarly, $q_{j} = 0$ for all $j \in J$. \qed
\enddemo
\proclaim{Claim 2} $S$ contains at most 19 distinct $\Bbb P^{1}$'s. 
\endproclaim 
\demo{Proof} Recall that $(H^{1,1}, (*.*))$ is 
of dimension 20 and of signature $(1, 19)$.                    
Assume that $S$ contains more than or equal to 20 distinct $\Bbb P^{1}$'s. 
Let $C_{1}, ..., C_{20}$ be 20 $\Bbb P^{1}$'s among them. 
Then, since $\Bbb R \langle [C_{1}], ..., [C_{20}] \rangle 
\subset NS(S) \otimes \Bbb R 
\subset H^{1,1}$ and $\text{dim}_{\Bbb R}\Bbb R \langle [C_{1}], ..., 
[C_{20}] \rangle = 20 = \text{dim}_{\Bbb R} H^{1,1}$ by Claim 1, 
we get $NS(S) \otimes \Bbb R = H^{1,1}$. However, 
$NS(S) \otimes \Bbb R$ is of signature $(0, 20)$ by (1) while 
$H^{1,1}$ is of signature $(1, 19)$, a contradiction. \qed 
\enddemo
Now we are done. \qed 
\enddemo
\demo{Remark} For each integer $m$ such that 
$0 \leq m \leq 19$, there actually exists a K3 surface of $a(S) = 0$ 
which contains exactly $m$ distinct $\Bbb P^{1}$'s and no other curves. 
\enddemo 
\demo{Construction} By [OZ], there exists a projective K3 surface 
$T$ which contains 19 $\Bbb P^{1}$'s, say, $C_{1}, ... , C_{19}$ 
whose intersection matrix $(C_{i}.C_{j})$ is of type $A_{19}$. 
Let $f : \Cal X \rightarrow \Cal B$ be the Kuranishi family of $T$ 
and identify the base space $\Cal B$ with an open set $\Cal U$ of 
the period domain $\Cal P$ of the K3 surfaces under some marking 
$\tau : R^{2}f_{*}\Bbb Z \simeq \Lambda_{K3} \times K$:
$$ \Cal B \simeq \Cal U \subset 
\Cal P := \{[\omega] \in \Bbb P(\Lambda_{K3} \otimes \Bbb C) 
\vert (\omega.\omega) = 0, \, (\omega.\overline{\omega}) > 0\}.$$
Let $c_{i}$ be the element of the K3 lattice $\Lambda_{K3}$ which 
corresponds to the class $[C_{i}]$ under the marking $\tau$. 
Define the subset $c_{i}^{\perp} \subset \Cal U$ by 
$c_{i}^{\perp} := \{[\omega] \in \Cal U \vert (\omega.c_{i}) = 0\}$. 
Let $0 \leq m \leq 19$ and choose a very general point $P$ 
of the space $c_{1}^{\perp} \cap ... \cap c_{m}^{\perp}$. 
(This space is of dimension $20 -m > 0$ by (2.4)). Then the fiber 
$\Cal X_{P}$ is a K3 surface which contains exactly $m$ distinct 
$\Bbb P^{1}$'s whose intersection matrix is of type $A_{m}$ 
and has no other curves. This also 
implies $a(\Cal X_{P}) = 0$. \qed \enddemo
\proclaim{Proposition (2.5)} Let $S$ be a complex torus of dimension 2. 
Assume that $a(S) = 0$. Then, $NS(S) \otimes \Bbb R$ is 
negative semi-definite with respect to $(*.*)$. 
\endproclaim 
\demo{Proof} Obvious. \qed \enddemo
 
\head 3. Proof of Huybrechts' Criterion \endhead 
\demo{Proof of the ``only if'' part} Any ample class gives a desired point. \qed 
\enddemo 
\demo{Proof of the ``if'' part} 
Let $S$ be a K\"ahler surface which has an inner integral point of 
$\Cal K^{*}(S)$. It is sufficient to show that $a(S) \not= 0, 1$. 
\enddemo

\proclaim{Lemma (3.1)} $a(S) \not= 1$.\endproclaim 
\demo{Proof} Assume to the contrary that $a(S) = 1$ and take the 
algebraic reduction $f : S \rightarrow C$, which, 
in the surface case, is a surjective morphism to a non-singular curve 
with connected fibers.
Let $F$ be a general fiber of $f$ and set $f(F) = P$.  
Then $[F] = [f^{*}(P)] \in NS(S)$, where $P$ is regarded as a divisor on $C$. 
Since $P$ is ample on $C$, 
the class $[P]$ is represented by a positive definite real $d-$closed 
$(1,1)$ form $\theta$. Set $\Theta := f^{*}\theta$.  
Then $[F] = [\Theta]$ and $\Theta$ is positive semi-definite 
at each point $Q \in S$. Therefore 
for a K\"ahler form $\omega$ and for any $\epsilon > 0$, 
we have $[\Theta + \epsilon \omega] \in \Cal K(S)$. Thus,  
$[F] = \lim_{\epsilon \rightarrow 0}[\Theta + \epsilon\omega] 
\in \overline{\Cal K(S)}$. Moreover, since 
$([F]. [\omega]) = \int_{F} \omega > 0$, we see that $[F] \not= 0$. 
Let $M$ be an inner integral point of $\Cal K^{*}(S)$. Then, by (1.6), 
we have $(M.[F]) > 0$, whence 
$(M + n[F])^{2} = M^{2}+ 2n(M.[F]) > 0$ for a large integer 
$n$. Since $M + n[F] \in NS(S)$, this implies $a(S) = 2$, 
a contradiction. \qed 
\enddemo 
The next Lemma reduces our problem to the case of minimal surfaces. 

\proclaim{Lemma (3.2)} Let $\tau : S \rightarrow T$ be the blow down 
of a $(-1)$-curve $E$. Then, 
\roster
\item $S$ is projective if and only if $T$ is projective.
\item $S$ is K\"ahler if and only if $T$ is K\"ahler. 
\item Assume that there exists an inner integral point 
$x$ of $\Cal K^{*}(S)$. Then there also exists an inner integral point 
of $\Cal K^{*}(T)$. 
\endroster
\endproclaim
\demo{Proof} The assertions (1) and (2) are well known. (However, it might 
be worth reminding here that the ``only if'' part of 
both (1) and (2) is false in general if dimension is three or higher 
and the center is not a point. 
One of instructive counterexamples is found in [Og].) 
\par 
Let us show the assertion (3). Recall that $H^{2}(S, K) =
\tau^{*}H^{2}(T,K) \oplus K[E] \simeq H^{2}(T,K) \oplus K[E]$ 
for $K = \Bbb Z, \Bbb R$. Moreover, this equality 
and isomorphism are compatible with the cup product and 
the Hodge decompositions. Let us regard $H^{1,1}(S)$ 
as a normed space by the product norm of $H^{1,1}(T)$ and 
$\Bbb R[E]$. Set $e := [E]$. Then the inner integral point $x \in 
\Cal K^{*}(S)$ is of the form $x = \tau^{*}y + ae$ where $y \in NS(T)$ and 
$a \in \Bbb Z$. We show that $y$ is an inner point of $\Cal K^{*}(T)$. 
Let $\sigma 
\in \Cal K(T)$. Then $\tau^{*}\sigma \not= 0$ and 
$\tau^{*}\sigma - \epsilon e \in \Cal K(S)$ for all sufficiently small 
positive 
real numbers $\epsilon$. Therefore $\tau^{*}\sigma \in \overline{\Cal K(S)}$. 
Since $x$ is an inner point of $\Cal K^{*}(S)$, 
there exists $r > 0$ such that $(x. \eta) \geq r \vert \eta \vert$ 
for all $\eta \in \overline{\Cal K(S)}$ by (1.6). In particular, 
$(x.\tau^{*}\sigma) \geq r \vert \tau^{*}\sigma\vert$. 
On the other hand, using $x = y + ae$ and applying the projection formula, 
we calculate 
$(x.\tau^{*}\sigma) = (y. \sigma)$. This 
together with the compatibility of the norms implies  
$(y. \sigma) \geq r \vert \sigma \vert$ for all $\sigma 
\in \Cal K(T)$, hence for all $\sigma \in \overline{\Cal  K(T)}$. \qed 
\enddemo
By virtue of (2.3), (3.1) and (3.2), in order to conclude 
the ``if'' part, it is now sufficient to show the 
following:
\proclaim{Lemma (3.3)} 
\roster
\item Let $S$ be a K3 surface. Assume that 
$\Cal K^{*}(S)$ contains an inner integral point $x$. 
Then $a(S) \not= 0$. 
\item Let $S$ be a complex torus of dimension 2. Assume that 
$\Cal K^{*}(S)$ contains an inner integral point $x$. 
Then $a(S) \not= 0$.
\endroster
\endproclaim
\demo{Proof of (1)} Assume to the contrary that $a(S) = 0$. 
Let $C_{1}, ... C_{m}$ ($0 \leq m \leq 19$) 
denote the distinct smooth rational curves on $S$ ((2.4)(2)). 
We argue dividing into two cases:
\par 
Case 1. $x \in \Bbb R \langle [C_{1}], ... , [C_{m}] \rangle$;
\par 
Case 2. $x \not\in \Bbb R \langle [C_{1}], ... , [C_{m}] \rangle$. 
\par
\vskip 4pt
Case 1. By (2.4), the subspace of $H^{1, 1}$
$$[C_{1}]^{\perp} \cap ... \cap [C_{m}]^{\perp}$$ 
is of signature $(1, 19 - m)$ (where $\perp$ is 
taken with respect to $(*.*)$). Therefore 
$$[C_{1}]^{\perp} \cap ... \cap [C_{m}]^{\perp} \cap \Cal C^{+}(S) \not= 
\emptyset.$$ 
Let $\eta$ be an element of this set. 
Then by (2.1), $\eta \in \overline{\Cal K(S)}$ and $\eta \not= 0$. 
On the other hand, by our assumption, we have 
$(x.\eta) = 0$. This contradicts (1.6). 
\par
\vskip 4pt
Case 2. In this case $m \leq 18$. Indeed, otherwise we would have
$\Bbb R \langle x, [C_{1}], ..., [C_{19}] \rangle = NS(S) \otimes \Bbb R = 
H^{1,1}$ and would get the same contradiction 
as in Claim 2 of (2.4). Therefore the subspace 
$$x^{\perp} \cap [C_{1}]^{\perp} \cap ... 
\cap [C_{m}]^{\perp}$$ 
is of signature $(1, 19 - m -1)$ and then
$$x^{\perp} \cap [C_{1}]^{\perp} \cap ... \cap [C_{m}]^{\perp} 
\cap \Cal C^{+}(S) \not= \emptyset.$$ 
Let $\eta$ be an element of 
$x^{\perp} \cap [C_{1}]^{\perp} \cap ... \cap 
[C_{m}]^{\perp} \cap \Cal C^{+}(S)$. Then by (2.1), $\eta \in 
\overline{\Cal K(S)}$ 
and $\eta \not= 0$. On the other hand, by the choice of $\eta$, we have 
$(x.\eta) = 0$, a contradiction. \qed \enddemo

\demo{Proof of (2)} Note that $(H^{1,1}, (*.*))$
is non-degenerate and of signature $(1, 3)$. 
Assume to the contrary that $a(S) = 0$. Then $x^{2} \leq 0$
and $x \not= 0$ by (2.5). We argue dividing into two cases: 
\par
Case 1. $x^{2} = 0$;
\par 
Case 2. $x^{2} < 0$.
\par 
\vskip 4pt
Case 1. Since $x^{2} = 0$ and $(x.\eta) > 0$ for all 
$\eta \in \Cal K(S)$, we have $x \in \overline{\Cal C^{+}(S)}$, 
whence $x \in \overline{\Cal K(S)}$ by (2.2). 
However, then $x^{2} = (x.x) > 0$ by (1.6), a contradiction. 
\par
\vskip 4pt
Case 2. Since $x^{2} < 0$, the subspace $x^{\perp} 
\subset H^{1,1}$ is of index $(1,2)$. Combining this with (2.2), 
we have
$x^{\perp} \cap \Cal K(S) = x^{\perp} \cap \Cal C^{+}(S) \not= \emptyset$.
Therefore there exists an element $\eta \in \Cal K(S)$ such that 
$(x. \eta) = 0$. However, this contradicts $x \in \Cal K^{*}(S)$. \qed
\enddemo

\vfill \eject 
\head References \endhead

\par
\vskip 20pt
Keiji Oguiso, Math. Inst. d. Univ. Essen, D-45117 Essen, Germany
\par 
E-mail: mat9g0\@spi.power.uni-essen.de 
\par
\vskip 10pt
Thomas Peternell, Math. Inst. d. Univ. Bayreuth, D-95445 Bayreuth, Germany
\par
E-mail: thomas.peternell\@uni-bayreuth.de  
\enddocument